\documentclass[12pt,titlepage]{article}
\usepackage{graphicx,color}
\usepackage{amsmath,amssymb,amsthm}
\input epsf

\newtheorem{thm}{Theorem}

\newtheorem{cor}{Corollary}

\newtheorem{prop}{Proposition}

\newtheorem{defi}{Definition}

\def\beq{\begin{equation}}\def\eeq{\end{equation}}
\def\beqn{\begin{eqnarray}}\def\eeqn{\end{eqnarray}}

\def\qed{\ifhmode\unskip\nobreak\fi\quad\ifmmode\Box\else$\Box$\fi}

\title{On colorful edge triples in edge-colored complete graphs}

\author{{\bf G\'abor Simonyi}$^{a,b,}$\thanks{Research partially supported by the
    National Research, Development and Innovation Office (NKFIH)
    grants K--116769 and K--120706. This work is also connected to the scientific program of the "Development of quality-oriented and harmonized R+D+I strategy and functional model at BME" project, supported by the New Hungary Development Plan (Project ID: TÁMOP-4.2.1/B-09/1/KMR-2010-0002).
}\\ \\
$^a$Alfr\'ed R\'enyi Institute of Mathematics, Hungarian Academy of
Sciences\\ \\
$^b$Department of Computer Science and Information Theory,\\ Budapest University of Technology and Economics\\ \\
 {\tt simonyi@renyi.hu}
}
\date{}

\begin{document}
\maketitle
\begin{abstract}
An edge-coloring of the complete graph $K_n$ we call $F$-caring if it leaves no
$F$-subgraph of $K_n$ monochromatic and at the same time every subset
of $|V(F)|$ vertices contains in it at least one completely multicolored version of $F$.
For the first two meaningful cases, when $F=K_{1,3}$ and $F=P_4$ we determine
for infinitely many $n$ the minimum number of colors needed for an $F$-caring
edge-coloring of $K_n$.
An explicit family of
$2\lceil\log_2 n\rceil$ $3$-edge-colorings of $K_n$ so that every quadruple of
its vertices contains a totally multicolored $P_4$ in at least one of them is
also presented.
Investigating related Ramsey-type problems we also
show that the Shannon (OR-)capacity of the Gr\"otzsch graph is strictly larger
than that of the five length cycle.

\bigskip
\bigskip
\par\noindent
%{\em Keywords:} Ramsey colorings, Kirkman triple systems, Shannon capacity, Mycielski construction

\end{abstract}

\section{Introduction}
\message{Introduction}

Erd\H{o}s \cite{Erdos} formulated the following Ramsey type problem: What is
the least number $f(n,p,q)$ of colors needed for an edge-coloring of the
complete graph
$K_n$ if it has the property that every subset of $p$ vertices spans at least
$q$ distinct colors? Initial investigations by Elekes, Erd\H{o}s, and F\"uredi
reported in \cite{Erdos} were followed by a more systematic study by Erd\H{o}s
and Gy\'arf\'as in \cite{ErdGyar}. A more general variant of the problem was
studied by Axenovich, F\"uredi, and Mubayi \cite{AxFurMub}, where they color
the edges of a not necessarily complete graph $G$ and require that every copy
of another graph $H\subseteq G$ receive at least a given number of colors.
Axenovich and Iverson \cite{AxIv} investigated the mixed anti-Ramsey numbers
that are the maximum and minimum numbers of colors to be used in an
edge-coloring of $K_n$ that avoids both monochromatic copies of a fixed graph
$G$ and also totally multicolored copies of another fixed graph $H$.

Here we initiate the investigation of yet another variant that will also be
connected to some of the above versions of the classical Ramsey problem. Let $g(n,F)$ denote
the minimum number of colors needed in an edge-coloring of $K_n$ if it
contains no monochromatic copy of $F$, furthermore, it contains at least one
totally multicolored copy of $F$ on every subset of $|V(F)|$
vertices.

If $F$ has only two edges then the two conditions coincide and we get
well-known numbers: for $F$ being two incident edges we need a proper edge
coloring of $K_n$, while in case of two non-incident edges we ask the
chromatic number of the complementary graph of the line graph of $K_n$, or
equivalently, the chromatic number of the special Kneser graph ${\rm KG(n,2)}$
(cf. \cite{LovKne}), which is $n-2$. The first meaningful cases of our problem
thus belong to graphs $F$ with $3$ edges and $4$ vertices,
since if $F$ is a triangle, or a complete graph, in general, then we are back
to the trivial problem of determining $f(n,p,{p\choose 2})={n\choose 2}$, with
$p\ge 3$, cf. \cite{ErdGyar}. Thus the first interesting cases belong to
$F=K_{1,3}$, and $F=P_4$, the star and the path on four vertices. We will
determine
$g(n,K_{1,3})$ and $g(n,P_4)$ for infinitely many $n$ in Section~\ref{sect:K13P4}.

Completely multicolored edge triples are investigated in a somewhat different
context in \cite{KS} and also in the subsequent paper \cite{BFGyR}. There we
have only three colors but
we can color in several rounds.  Vera T. S\'os asked in 1991, how many
$3$-edge-colorings of $K_n$ ($n\ge 3$) are needed to have every triangle $3$-colored in
at least one of them. This question was related to a special case of the
so-called perfect hashing problem, cf. \cite{FK,KHGP,KM}.
In \cite{KS} we gave a construction proving that the requested
minimum number is at most $\lceil\log_2n\rceil-1$, while an easy argument
shows that it is bounded from below by $\lceil\log_3(n-1)\rceil$.
Different edge triples were also considered in \cite{KS}, but we did not find meaningful bounds in the case when these edges form a $P_4$. In Section~\ref{sect:rounds} we give upper bounds for the number of $3$-edge-colorings needed if we want that at least one $P_4$ gets $3$-colored on every $4$ vertices in at least one of our colorings.
\par\noindent
Note that if we wanted to only $2$-color every triangle in several
rounds of $3$-edge-colorings then we would need exactly
$\lceil\log_3 k(n)\rceil$ rounds, where $k(n)=\min\{k: R(3;k)>n\}$, where
$R(3;k)$ is the minimum size of a complete graph that cannot be $k$-edge colored without creating a monochromatic triangle, i.e., the Ramsey number $R(3,3,\dots,3)$ with the
number of $3$'s being $k$. (We simply have to encode
each color with a ternary sequence to convert the coloring into $\lceil\log_3 k(n)\rceil$ $3$-edge-colorings.) The growth rate of this number is
well-known to be closely related to the Shannon (OR-)capacity of triangle-free
graphs. Elaborating on this connection we observe, using a result of Cropper,
Gy\'arf\'as, and Lehel \cite{CGyL}, that the answer to a famous question about
the growth
rate of $R(3:k)$ depends on whether the Shannon (OR-)capacity of Mycielski
graphs go to infinity. As a first modest step towards such investigations we
show that the Shannon (OR-)capacity of the Gr\"otzsch graph is strictly more
than that of $C_5$, the cycle of length five.  These questions are discussed
in the last subsection of Section~\ref{sect:rounds}.

\section{$K_{1,3}$- and $P_4$-caring colorings} \label{sect:K13P4}

Along with introducing some other related parameters, let us repeat the
definition of $g(n,F)$ that we already defined less formally in the
Introduction.

\begin{defi}\label{defi:szamok}
Let $F$ be any fixed graph and $n\ge |V(F)|$.
An edge-coloring of $K_n$ in which no copy of $F$ is monochromatic and at the
same time every subset of $|V(F)|$ vertices contains a totally multicolored
(rainbow) copy of $F$ is called $F$-caring.
Let $g(n,F)$ denote the minimum number of colors needed to an $F$-caring
coloring of $K_n$.

Let $b(n,F)$ denote the minimum number of colors needed for an edge-coloring
of $K_n$ that contains no monochromatic copy of $F$.

Finally, let $a(n,F)$ denote the minimum number of colors needed for an
edge-coloring
of $K_n$ that makes at least one copy of $F$ totally multicolored in
every subset of $|V(F)|$ vertices.
\end{defi}

Thus $b(n,F)$ is simply the smallest $k$ for which the Ramsey number
$R_k(F,F,\dots,F)>n$, where the index $k$ refers to the number of colors
used. It is immediate from the above definition that $g(n,F)\ge
\max\{a(n,F),b(n,F)\}$.
It is a trivial matter to determine $b(n,K_{1,3})$ exactly,
and it is also easy to tell the value of $b(n,P_4)$ for
infinitely many values of $n$.

\begin{prop}\label{prop:bicsill}
One has
$b(n,K_{1,3})=\left\lceil\frac{n-1}{2}\right\rceil$ for every $n\ge 4$.
\end{prop}

\proof
Any color class that contains no $K_{1,3}$
must give a graph with maximum degree at most $2$, thus it is a collection of
vertex-disjoint paths and cycles. It follows that the number of edges in each
color class is at most $n$ and thus the lower bound $b(n,K_{1,3})\ge
\lceil\frac{n-1}{2}\rceil$ follows.

On the other hand, it is well-known that if $n$ is odd, then the edge set of
$K_n$ can be partitioned into $\frac{n-1}{2}$ Hamiltonian cycles, see e.g.
\cite{Bondy}. Letting the edges of these cycles form the color classes the
Proposition follows for odd values of $n$.
%Note that if $n\equiv 3\pmod 6$ then we can also color
%with the parallel classes of Kirkman triple systems as in the proof of
%Proposition~/ref{prop:bip4}.

If $n$ is even, then we can simply take an optimal edge-coloring of $K_n$ with
$n-1$ colors, and pair up the colors in an arbitrary way. For each pair of
colors define a new color class, obviously, no vertex will connect to three
edges in it. This proves the required upper bound for even values of $n$.
\hfill$\Box$

\medskip
\par\noindent
For the next proposition we need some notions from the theory of block
designs. A Steiner triple system (STS for short) is a family of $3$-element
subsets of
a basic set of $n$ elements having the property that every pair of elements
belongs to exactly one triple of the family. An STS is called resolvable if
it can be partitioned into subfamilies, called {\em parallel classes}, in which each
element of the
basic set belongs to exactly one triple. (That is, in each parallel class the
triples are pairwise disjoint and their union is the entire basic set.)
Resolvable Steiner triple systems are called Kirkman triple systems (KTS for
short) referring to Kirkman's puzzle in the middle of the 19th century asking
for the
existence of such a system of triples and their parallel classes for $n=15$,
cf. \cite{Stinson}. In full
generality it was proven only more than a hundred years later by Ray-Chaudhury
and Wilson \cite{RCW} that Kirkman triple systems exist for all
$n\equiv 3 \pmod 6$. (It is easy to see that they cannot exist for any other
$n$.)

\begin{prop} \label{prop:bip4}
If $n\equiv 3\pmod 6$, then $b(n,P_4)=\frac{n-1}{2}.$
\end{prop}

\proof
Consider an edge-coloring of $K_n$ with $b(n,P_4)$ colors in which no $P_4$ is
monochromatic. Then each color class is a $P_4$-free graph, so it must
be a collection of vertex-disjoint triangles and stars. This implies that the
number of edges in each color class is at most $n$. So the number of color
classes $b(n,P_4)\ge \frac{n-1}{2}$.

For equality we need a coloring in which all color classes contain only
triangles (none of them are stars). As every edge gets exactly one color,
this requires that the union of the triangles forming the color classes be a
Steiner triple system if we consider the triangles as triples. Furthermore, as
each color class should consist of vertex-disjoint triangles, they should
form the parallel classes of a Kirkman triple system. Conversely, if a Kirkman
triple system on $n$ points exists, then we can use its parallel classes to
define color classes of the edges so that each edge is colored with some
identifier of the parallel class it belongs to (that is, which contains a
triple, that covers this edge). Such a coloring uses
$\frac{n-1}{2}$ colors and it lets no $P_4$ remaining monochromatic. This proves $b(n,P_4)\le\frac{n-1}{2}$ whenever $n\equiv 3 \pmod 6$.
\hfill$\Box$

\medskip
\par\noindent
The two propositions above give obvious lower bounds for the value of
$g(n,K_{1,3})$ and $g(n,P_4)$. Below we show that these lower bounds
are tight.

\begin{prop}\label{prop:gcsill}
If $n\equiv 3\pmod 6$, then
$$g(n,K_{1,3})=\frac{n-1}{2}$$.
\end{prop}

\proof
Since $g(n,F)\ge b(n,F)$ by definition, $g(n,K_{1,3})\ge\frac{n-1}{2}$ is
immediate from Proposition~\ref{prop:bicsill}.
For the reverse inequality take a coloring based on a Kirkman triple system as
in the proof of Proposition~\ref{prop:bip4}, for the numbers in the statement
such a system exists by the result of Ray-Chaudhury-Wilson \cite{RCW}. We show
this is a $K_{1,3}$-caring coloring.

It is clear that no $K_{1,3}$ is monochromatic, since each vertex has degree
$2$ only within each graph formed by edges of a color class. Note, that if two
incident edges belong to the same color class, then
the edge extending this pair of edges to a triangle should also belong to the same
color class. Now consider a set of four points $\{a,b,c,d\}$. If there is a monochromatic
triangle on them formed by, say, vertices $a,b,c$, then the fourth vertex must
be connected to each of $a,b$, and $c$ in a different color. (If two of them
were the same, then it would require an edge of the monochromatic triangle
to belong to another monochromatic triangle, which is impossible.) Thus the
$K_{1,3}$ spanned by edges $da,db,dc$ is rainbow. If there is no monochromatic
triangle on our $4$ points, then it cannot contain even two incident edges
with the same color and thus all the four $K_{1,3}$'s it contains are rainbow.
\hfill$\Box$

\medskip
\par\noindent
In what follows we will use the notation $a\oplus_nb$ for the modulo $n$ sum of integers $a$ and $b$.

\begin{thm}
If $n$ is a power of $3$, then $$g(n,P_4)=\frac{n-1}{2}$$.
\end{thm}

\proof
Since $g(n,F)\ge b(n,F)$, we have $g(n,P_4)\ge\frac{n-1}{2}$ from
Proposition~\ref{prop:bip4}. Now let $n=3^t$. For the reverse inequality we
show a Kirkman triple system that has the property that every quadruple
of vertices contains some $P_4$ the three edges of which belong to three different parallel classes of the KTS. (In what follows we will often refer to the parallel classes as colors. Colorings defined by the parallel classes of a KTS as color classes will be called KTS colorings. Note that no $P_4$ can be monochromatic in a KTS coloring according to the argument in the proof of Proposition~\ref{prop:bip4}).

It is straightforward to check that if a $K_4$ is edge-colored with at least $4$ colors then it must contain a rainbow $P_4$. If a $K_n$ is edge-colored with the parallel classes of a KTS then there is only one way not to have $4$ differently colored edges on a fixed set of four vertices, namely when it is colored with three colors and the coloring partitions its edge set into three matchings, each containing exactly two edges. (In other words it is an optimal proper edge coloring of this $K_4$.) This simply follows from the fact that if two incident edges have the same color in a KTS coloring, then the third edge with which they form a triangle must also belong to the same color class. So it is enough to show a KTS on $n=3^t$ vertices for every $t$ with the property that the KTS coloring it defines will not partition the edge set of any $K_4$ into three perfect matchings. Let us call a KTS with this property that it has at least $4$ parallel classes intersecting the edge set of any $K_4\subseteq K_n$ a good KTS.

The smallest non-trivial KTS that exists is on $9$ points and is easily shown
to have this property as follows. Assume that
for some four points $\{a,b,c,d\}$ the coloring it defines colors $\{a,b\}$
and $\{c,d\}$ the same color. Then the
edges $\{a,b\}$ and $\{c,d\}$ belong to disjoint triples
of the same
parallel class of the KTS. This needs two distinct vertices belonging to those
two triples outside the set $\{a,b,c,d\}$. The same can be said about the
other two pairs of edges, so altogether we need six vertices outside
$\{a,b,c,d\}$, but we only have $5$ as there are only $9$ points altogether. So the KTS on $9$ points is a good KTS.

Now we prove the theorem by induction on $t$. Assume we have a good KTS on $3^{t-1}$ points, and show that one exists also on $3^t$ points.
(The base case was just verified for $t=1$. In fact, we could have taken $t=0$ for the base case as
well, and get the KTS on $9$ points also by the induction step.)
Partition our $3^t$ points into three equal classes $A_j$ and label them as
$(\ell,j)$ with $\ell\in\{0,1,2\}$ referring to their partition class and
$j\in \{0,1,\dots, 3^{t-1}-1\}$ their position according to an arbitrary
labelling within their class.
Take three copies of the existing good KTS on $3^{t-1}$ points with parallel
classes $C_1, \dots, C_m, m=\frac{3^{t-1}-1}{2}$ on the three partition
classes and consider the union of their $i$th parallel classes as the $i$th
parallel class of the KTS we are constructing on $3^t$ points. For
$j=0,1,\dots,3^{t-1}-1$ we define
$C_{m+j+1}:=\{\{(0,a),(1,a\oplus_{3^{t-1}}j),(2,a\oplus_{3^{t-1}}2j)\}: 0\le
a<3^{t-1}\}$. Then $\{C_1,\dots,C_{\frac{3^t-1}{2}}\}$ is a Kirkman triple system, we show it is good.

Call a $K_4$ bad (with respect to a KTS) if its edges belong to only $3$
parallel classes of the KTS. Consider any four points, we want to show that
they cannot induce a bad $K_4$ (with respect to the KTS just defined). If the
four points are all in the same class $A_j$ then they could induce a bad $K_4$
only if the KTS on $3^{t-1}$ points we started with would not be a good one,
which is not the case. If our four points are in $3$ different classes, then
exactly one pair of these four points belongs to the same class and thus to a
$C_k$ with $0\le k\le m$ with $m$ still meaning the value
$\frac{3^{t-1}-1}{2}$. Then this pair of points forms a differently colored
edge then the complementary two points in our $4$-element set, so the edge set of our $K_4$
cannot be partitioned into three perfect matchings and thus this $K_4$ cannot be bad. If the four points are in two classes so, that one is in one
class, say $A_j$, and three of them are in the other class, say $A_h$ then the
three edges of our $K_4$ that has one endpoint in $A_j$ and the other in $A_h$
all belong to different $C_k$'s with $m<k$ while the remaining three edges all
belong to some $C_k$'s with $k\le m$. So again, we have at least $4$ colors in
our $K_4$. Lastly, if two of our points are in $A_j$ and two are in $A_h$ with
$j\ne h$, then the four edges running between the two classes must belong to
at least three different $C_k$'s with $k>m$. This is because only pairs of independent edges among them can have the same color and if $a,b\in A_j$, $c,d\in A_h$, and edges $\{a,c\}, \{b,d\}$ have the same color, then we
have $c-a\equiv d-b\ (\textrm{mod}\ 3^{t-1})$. Both $\{a,d\}$ and $\{b,c\}$
should then belong to a different color class, and these two may not be equal,
otherwise we also have  $c-b\equiv d-a\ (\textrm{mod}\ 3^{t-1})$, and the last
two congruences would imply $b-a\equiv a-b\ (\textrm{mod}\ 3^{t-1})$
contradicting that $a$ and $b$ are different.
The remaining two edges $\{a,b\}$ and $\{c,d\}$ belong to $C_k$'s with
$k\le m$, so we again have at least four differently colored edges in our $K_4$.
This completes the proof of the theorem.
\hfill$\Box$

\bigskip
\medskip

\medskip
\par\noindent
{\it Remark 1.} We have $g(n,F)\ge\max\{a(n,F),b(n,F)\}$ for every graph and we have seen
above that for $K_{1,3}$ and $P_4$ we have $g(n,K_{1,3})=b(n,K_{1,3})$ and
$g(n,P_4)=b(n,P_4)$. In fact, $a(n,K_{1,3})$ and $a(n,P_4)$ are indeed
strictly smaller than $b(n,K_{1,3})$ and $b(n,P_4)$. This can be seen by
using that $a(n,K_{1,3}),a(n,P_4)\le f(n,4,4)$ (for the meaning of $f(n,p,q)$ see the beginning of the Introduction), and Erd\H{o}s and
Gy\'arf\'as's result from \cite{ErdGyar} stating that $f(n,4,4)\le O(n^{2/3})$,
improved by Mubayi \cite{Mubayi} (via an explicit construction) to
$f(n,4,4)\le n^{1/2}e^{c\sqrt{\log n}}$ with an absolute constant $c>0$.
Nevertheless, $a(n,F)>b(n,F)$ can also happen. This is trivially so if $F$ is
a complete graph, but also for any $F$ with $|E(F)|>\frac{1}{2}{|V(F)|\choose 2}$
we immediately have $a(n,F)\ge b(n,F)$. For $F$ being a complete graph minus
one edge we also trivially have
$g(n,f)=a(n,F)=f(n,|V(F)|,{|V(F)|\choose 2}-1)$. In particular,
$f(n,4,5)\ge\frac{5}{6}(n-1)$ is proven in \cite{ErdGyar} implying that for
$F=K_4\setminus\{e\}$ (where $e$ is an edge of $K_4$) we have
$b(n,F)\le f(n,4,3)<f(n,4,5)=a(n,F)$. \hfill$\Diamond$

\smallskip
\par\noindent
{\it Remark 2.} Although an analogous result to that of Ray-Chaudhury and Wilson \cite{RCW}
about the existence of Kirkman triple systems
was proven by Hanani, Ray-Chaudhury, and Wilson \cite{HRCW} for quadraple
systems (that is the existence of resolvable so-called balanced incomplete
block designs of block size $4$ whenever a trivial necessary condition is met),
it cannot be directly used for generalizing even Proposition~\ref{prop:gcsill}
to determine $g(n,K_{1,4})$. This is because a colouring based on such
designs will not guarantee the $4$-coloring of a $K_{1,4}$ subgraph on every
$5$ vertices.
\hfill$\Diamond$

\section{Coloring in several rounds}
\label{sect:rounds}

\subsection{Relaxations of perfect hashing}

The perfect hashing problem, investigated in \cite{FK, KHGP, KM} asks for the minimum length of sequences over an alphabet of size $b\ge k$, such that one can give $n$ such sequences so that for any $k$ of them there is a position where all of these $k$ sequences differ. While the answer is trivially $\log_bn$ if $k=2$, for $k=3$ the answer is not known even for the smallest meaningful alphabet size $b=k=3$. Denoting the required minimum length by $t(n,b,k)$ in general, it is known that $t(n,b,k)=\Theta(\log n)$ for any fixed pair $(b,k)$ and the relevant quantity asked for is $c_{b,k}=\liminf_{n\to\infty}\frac{t(n,b,k)}{\log n}$. The best known bounds on $c_{b,k}$ in general are given in \cite{KM}, in case of $c_{3,3}$ their values are

$$\frac{1}{\log\frac{3}{2}}\le c_{3,3}\le\frac{4}{\log\frac{9}{5}},$$
where, as always in the sequel when not stated otherwise, the logarithms are
meant to be on base $2$. Numerically the above bounds mean $1.709\le c_{3,3}\le 4.717.$

We remark that the upper bound is implicit in Elias \cite{PE}, and that for $c_{4,4}$ some improvements were found recently by Dalai, Guruswami, and Radhakrishnan \cite{Radha}. From now on we are focusing here only on the $b=k=3$ case.

While this problem is notoriously difficult, it is interesting to look at
relaxations where we do not require all triples of our sequences being
pairwise different at some position (this is the relation we called
"trifference" in \cite{KS}) but only a selected collection of all triples
should have this property. Such a selected collection is most naturally defined by some underlying
structure. This explains what could be the motivation for Vera T. S\'os in
1991 for asking how long our sequences should be if the basic set is identified to
the edge set of $K_n$ and only those triples should be trifferent at some
position that form a triangle. In other words: How many $3$-edge-colorings of
$K_n$ are needed to make every triangle $3$-colored in at least one of them?
Though we still do not know the exact answer to this question, it indeed
seems to be more easily tractable than the problem of $c_{3,3}$, because here
we could at least give an explicit construction (see \cite{KS}) that implies a
$\lceil\log_2n\rceil-1$ upper bound on the length of sequences needed, while
an obvious lower bound is $\lceil\log_3(n-1)\rceil$. These two bounds
differ first at $n=9$ for which Blokhuis, Faudree, Gy\'arf\'as, and Ruszink\'o \cite{BFGyR} proved that the upper bound is tight. For a completely
solved variant involving orientations in place of colorings, see \cite{HS}. (For the general
trifference problem the best known upper bounds on $c_{3,3}$ are obtained via
the probabilistic method, that is, they are not constructive.) If
instead of the edges of triangles we want the $3$ edges of every $K_{1,3}$ or
every $3$ independent edges to become trifferent in at least one coloring,
then the required minimum number of colorings will be $c_{3,3}\log n+o(\log
n)$, see \cite{KS}. (The reason is simple: In case of $K_{1,3}$, we have to
make all the $(n-1)$ edges incident to a fixed vertex have the property that
any three of them is trifferent at some coloring, while encoding the $(n-1)$
or $n$ colors of an optimal proper edge-coloring this way will do job; and a
similar reasoning works for three independent edges.)

We did not find non-trivial bounds in \cite{KS} for the triples given by the
edge sets of $P_4$'s. A trivial upper bound is $2c_{3,3}\log n+o(\log n)$,
since with that many colorings we can make all the ${n\choose 2}$ edges having
the property that any triple is trifferent at some coloring. We also have
$\log_3(n-1)$ as a trivial lower bound as any two edges incident to the same
fixed vertex should get a different color in some coloring as they appear
together in some $P_4$. (Note that in this kind of problems it is typically
the constant multiplier of $\log n$ that we are seeking for.)

A consequence of Erd\H{os} and Gy\'arf\'as's \cite{ErdGyar} and of Mubayi's result \cite{Mubayi} already mentioned in the first remark at the end of the previous section is that if we require only that the three edges of at least one $P_4$ of every $K_4$ subgraph get trifferent in some $3$-edge-coloring then the number of colorings we need is strictly less than the above mentioned trivial upper bound.

\begin{defi}\label{defi:pn}
Let $p(n)$ denote the minimum number of $3$-edge-colorings of $K_n$ needed if
we want that for all $4$ vertices the edge set of some $P_4$ subgraph on those
$4$ vertices gets totally multicolored in at least one of these colorings.
\end{defi}

\begin{prop}\label{prop:egyP4}
$$p(n)\le \frac{1}{2}c_{3,3}\log n+o(\log n).$$
\end{prop}

\proof
Color the edges of $K_n$ with $m:=n^{1/2}e^{c\sqrt{\log n}}=n^{1/2+o(1)}$ colors with Mubayi's method so that at least $4$ colors appear on every $K_4$ subgraph. This means that at least one $P_4$ of each $K_4$ is totally multicolored. Now encode each color with a ternary sequence of length $t(m,3,3)$ so that any $3$ sequences are trifferent in at least one coordinate. Now substituting every color in the previous coloring with a sequence of $3$-colorings according to the ternary codes of colors, we will get a required sequence of $3$-colorings, since any three edges that received three different colors in the first coloring must get three different colors in the $3$-coloring that belongs to the position where the codewords of those $3$ colors are trifferent.

Since $t(m,3,3)=c_{3,3}\log m+o(m)=c_{3,3}\frac{1}{2}\log n+o(\log n)$ the result follows.
\hfill$\Box$

\medskip
\par\noindent
{\it Remark 3.} A $K_4$ that is edge-colored with $4$ colors contains at least $4$ totally multicolored $P_4$'s out of the altogether $12$ ones. It would be interesting to see whether one could apply Mubayi's coloring $3$ times in such a way that all $P_4$'s get $3$-colored in at least one of them. If this would be possible, then the argument in the proof of Proposition~\ref{prop:egyP4} would give that with $\frac{3}{2}c_{3,3}\log n+o(\log n)$ $3$-edge colorings of $K_n$ one can totally multicolor all its $P_4$'s thus providing a strictly better upper bound than the trivial $2c_{3,3}\log n+o(\log n).$
\hfill$\Diamond$

\medskip
In the following theorem we use a variation of an idea from \cite{KS} to
give an explicit construction that provides another upper bound for $p(n)$. Because
of the lack of our knowledge of the value of $c_{3,3}$, we do not know whether
this upper bound is better or worse than the one provided by the previous
proposition.

\begin{thm}\label{thm:binnegy}
$$p(n)\le 2\lceil\log_2 n\rceil.$$
\end{thm}

\proof
First we give $\lceil\log_2 n\rceil$ $4$-colorings of the edges of $K_n$ with
the property that any $4$-element subset of the vertex set contains a $P_4$
subgraph that is totally multicolored in at least one of them.

Label the $n$ vertices with the binary form of the numbers $0,1,\dots, n-1$
padded with $0$'s at the beginning if necessary to make all these binary
sequences having the same length $t:=\lceil\log_2 n\rceil$.
This way each vertex is attached a different $0-1$ sequence of length
$t$.

\par\noindent
Let $u,v$ be two distinct vertices and $u_i,v_i$ be the $i$th coordinate in
the binary sequences attached to $u$ and $v$, respectively. We define a coloring
of edge $u,v$ for every $1\le i\le t$ as follows.
Color the edge $\{u,v\}$

blue if $u_i=v_i=0,$

green if $u_i=v_i=1,$

red\ if $u_i\neq v_i$ and $\forall j<i: u_j=v_j,$

yellow if $u_i\neq v_i$ and $\exists j<i: u_j\neq v_j.$

\smallskip
\par\noindent
Consider any set of $4$ vertices $a,b,c,d$. If there is any coordinate $i$
where exactly two of $a_i, b_i, c_i, d_i$ are $0$ and two are $1$, then we have
two non-incident edges one of which is blue, while the other is green, and the
remaining $4$ edges of the $K_4$ with vertices $a,b,c,d$ are all red and
yellow. So those $P_4$'s (there are $4$ of them), that contain the blue and
green edges are totally multicolored.

If there is no such coordinate where exactly two of these four sequences are
$0$, then there are at least three such coordinates where one of the values is
different from the other three. Without loss of generality we may assume that
$a_i\neq b_i=c_i=d_i, b_j\neq a_j=c_j=d_j$, and $c_h\neq a_h=b_h=d_h$, where
we may also assume $i<j<h$ and that these $i,j,h$ are the smallest coordinates
where the above relations hold. Then the coloring rule gives that in the $j$th
coloring the edge $\{a,b\}$ is yellow, $\{b,c\}$ is red, and $\{c,d\}$ is
either blue or green, thus giving a totally multicolored $P_4$. In fact, also
the $P_4$ with edges $\{a,b\}, \{b,d\}, \{c,d\}$ is totally multicolored in
the $j$th coloring, while two other $P_4$'s will be totally multicolored in
the $h$th coloring, so altogether we can guarantee at least $4$ different
totally multicolored $P_4$'s on any given quadruple of vertices.

Since $t(4,3,3)=2$ (just take the $2$-length sequences $00, 01, 12, 22$, any
three of these are trifferent in at least one of the two coordinates), we can encode
the $4$ colors with two length ternary sequences that make any three of them
trifferent in at least one coordinate. Thus replacing each of the
$4$-colorings in the previous construction with two appropriate $3$-colorings
(just identifying the colors blue and green in one and the colors yellow and
red in the other), we obtain $2\lceil\log_2 n\rceil$ $3$-edge colorings
satisfying that for every $4$ vertices some $P_4$ on them is totally
multicolored by at least one of these colorings. In fact, we have seen,
that at least $4$ out of the $12$ $P_4$'s on any given four vertices will be rainbow at some coloring.
\hfill$\Box$

\medskip
\par\noindent
Naturally, a similar remark to that of Remark 3. above could be formulated
concerning the coloring used in the proof of Theorem~\ref{thm:binnegy}.

\medskip
\par\noindent
It should be clear from the proof of Theorem~\ref{thm:binnegy} that if we label the vertices of $K_n$ with $t$-length binary sequences so, that for any four of them there is a coordinate where exactly two of these four sequences contain a $1$, then $p(n)\le t$. (We could simply use the same coloring strategy but identifying colors red and yellow.) Alon, K\"orner, and Monti \cite{AKM} investigated the maximum number of binary sequences of some given length with this property. Their results give, that if $\mu(n)$ is the minimum $t$ for which $n$ such $t$-length sequences can be given, then $\frac{1}{0.78}\log_2 n+o(\log n)\le\mu(n)\le \frac{3}{\log_2\frac{8}{5}}\log_2 n+o(\log n)\approx 4.425\log_2 n+o(\log n)$. In a similar manner one can obtain an upper bound on the minimum number of $3$-colorings needed to make every $P_4$ totally multicolored in at least one of them if we use so-called $(2,2)$-separated binary sequences that means a set of sequences with the property that for any two disjoint pair of them, there is a coordinate where both sequences of the first pair contain a $0$ and both sequences of the other pair contain a $1$. The problem of maximum cardinality of such a set of binary sequences of some fixed length is considered, for example, in \cite{FK} and \cite{KS1}.

\smallskip
\par\noindent
An immediate consequence of Proposition~\ref{prop:egyP4},
Theorem~\ref{thm:binnegy}, and the discussion in the previous paragraph is the following. Since we do not know whether
$c_{3,3}$ is smaller or larger than $4$ (cf. the bounds quoted at the beginning of this section for $c_{3,3}$), and the value of $\mu(n)$ is also unknown, we do not know which expression
achieves the minimum in the statement.

\begin{cor}\label{cor;combined}
$$p(n)\le\min\{2\lceil\log_2 n\rceil, \frac{1}{2}c_{3,3}\log n+o(\log n), \mu(n)\}.$$
\hfill$\Box$
\end{cor}

\medskip
\par\noindent
It is rather annoying that we do not even know whether $p(n)=\Omega(\log n)$
is true. Since $p(n)\ge\log_3 a(n,P_4)$ and $a(n,P_4)\ge f(n,4,3)$, this would
immediately follow if we knew a lower bound $f(n,4,3)\ge n^{\varepsilon}$ for
any constant $\varepsilon>0$. This is, however, not known, in fact, the
behavior of $f(n,4,3)$ is mentioned as a most annoying open problem in
\cite{ErdGyar}, where a $c\sqrt{n}$ lower bound is given for it, while the
authors also mention that they could not even prove that $\frac{f(n,4,3)}{\log n}$ tends to infinity.

\subsection{Bicoloring the triangle and Shannon capacity} ~\label{subsect:Shannon}

It is rather trivial to see that if we want to two-color only certain subgraphs in several rounds then the problem of determining the minimum number of rounds is equivalent to some of those problems discussed in Section~\ref{sect:K13P4}, in particular in determining the value $b(n,F)$. Indeed, if we can use a fixed number $b$ of colors in several rounds and we want that each subgraph $F$ gets at least two colors in at least one of the rounds, then the minimum number of rounds needed is exactly $\lceil\log_b b(n,F)\rceil$.

Propositions~\ref{prop:bicsill}~and~\ref{prop:bip4} show that in case of
$F=K_{1,3}$ and $F=P_4$ we have $\log b(n,F)=\Theta(\log n)$. The situation is
rather different, however, for $F=K_3$, because $b(n,K_3)=\min\{k: R(3;k)>n\}$
and the growth rate of $R(3;k)$ is at least exponential in $k$. This means
that the minimum number of colorings (even with two colors) making every
triangle bicolored in at least one of them is bounded from above by
$O(\log\log n)$. On the other hand, it is not known, whether the growth rate
of $R(3;k)$ is exponential or even faster in $k$. This is a famous open
problem of Erd\H{o}s, who, according to \cite{CG} (see Chapter 2.5) offered \$250 for determining the limit $\lim_{k\to\infty}R(3;k)^{1/k}$ and \$100 for deciding whether it is finite. (The best lower bound on this limit we know about is $(321)^{1/5}\approx 3.171765\dots$ due to Exoo \cite{Exoo}.)
This question is well-known to be equivalent to a question about Shannon capacity, cf. \cite{AO, EMcT, NR}.

\begin{defi}\label{defi:ORprod}
Let $G$ and $H$ be two graphs. Their OR-product $G\otimes H$ is defined as follows.
$$V(G\otimes H)=V(G)\times V(H)$$ and $$E(G\otimes H)=\{\{(g_1,h_1).(g_2,h_2)\}: (g_1,g_2)\in E(G)\ {\rm or}\ (h_1,h_2)\in E(H)\}.$$
$G^t$ denotes the $t$-fold OR-product of graph $G$ with itself.
\end{defi}

\begin{defi}\label{defi:ShORcap}
The Shannon OR-capacity of graph $G$ is defined as the always existing limit
$$C_{\rm OR}(G)=\lim_{t\to\infty}\sqrt[t]{\omega(G^t)},$$
where $\omega$ stands for the clique number.
\end{defi}

Shannon capacity is most often defined by another graph product (the so called
AND-product) and via independence numbers, see \cite{Sha, LovShan}, but the
two notions are essentially the same as the Shannon OR-capacity of a graph $G$
we defined here is simply the more usual Shannon capacity (the ``Shannon
AND-capacity'') of the complementary graph $\bar G$. We use this less often used version because certain statements are more natural in this language and it is more convenient for our goals.
(This approach also appears in several papers, cf. for example \cite{GKV, SaSi}.)
It is well-known that $\omega(G)\le C_{\rm OR}(G)\le\chi(G)$, where $\chi(G)$
is the chromatic number of $G$. In general, Shannon (OR)-capacity is very hard
to determine, the only odd cycle longer than $3$ for which it is known is
$C_5$, and the determination of $C_{\rm OR}(C_5)$ was a major result of
Lov\'asz \cite{LovShan}. It was also a sort of breakthrough when Bohman and
Holzman \cite{BH} proved with an ingenious construction that the Shannon
OR-capacity of every odd cycle is strictly larger than its trivial lower bound
$2$, i.e., $C_{\rm OR}(C_{2k+1})>2$ for every positive integer $k$. Erd\H{o}s,
McEliece, and Taylor \cite{EMcT} (and independently but later Alon and
Orlitsky \cite{AO}, cf. also \cite{NR}) showed, that the supremum of the
possible Shannon OR-capacities a triangle-free graph can have is equal to
$\lim_{k\to\infty}\sqrt[k]{R(3;k)}$. That is, Erd\H{o}s's question whether
this limit is infinite is equivalent to the question, whether a triangle-free
graph can have arbitrarily large Shannon OR-capacity. Based on a nice
observation by Cropper, Gy\'arf\'as, and Lehel \cite{CGyL} we claim that to
decide this, it would be enough to know the Shannon OR-capacity of the
well-known Mycielski graphs \cite{Myc} that we define below.

\begin{defi}\label{defi:Myc}
For any graph $G$ on vertex set $\{v_1,\dots,v_n\}$ its Mycielskian $M(G)$ is defined as follows.
$$V(M(G))=V(G)\cup\{u_i: v_i\in V(G)\}\cup\{z\},$$
$$E(M(G))=E(G)\cup\{\{u_i,v_j\}: \{v_i,v_j\}\in E(G)\}\cup\{\{u_i,z\}: v_i\in V(G)\}.$$
By the $k$th Mycielski graph $M_k$ we mean the result of the $(k-2)$th iteration of the operation $M(.)$ defined above when starting with $M_2:=K_2$. The vertex $u_i$ we refer to as the twin vertex of $v_i$.
\end{defi}

We index $M_k$ as defined here, because as is well-known, the Mycielski construction increases the chromatic number by $1$ (while not changing the clique number) and this way $M_k$ denotes the $k$-chromatic Mycielski graph. Note that $M_3$ is just $C_5$ and $M_4$ is the graph also called Gr\"otzsch graph.

\begin{thm}\label{prop:ShaMyc}
$$\lim_{k\to\infty}\sqrt[k]{R(3;k)}=\lim_{r\to\infty}C_{\rm OR}(M_r).$$
\end{thm}

\proof
By the result of Erd\H{o}s-McEliece-Taylor \cite{EMcT} and Alon-Orlitsky \cite{AO} already quoted above we know that $$\lim_{k\to\infty}\sqrt[k]{R(3;k)}=\sup\{C_{\rm OR}(G): K_3\nsubseteq G\}.$$

\par\noindent
Cropper, Gy\'arf\'as, and Lehel \cite{CGyL} observed that every triangle-free graph is an induced subgraph of $M_r$ if $r$ is large enough. Since $C_{\rm OR}(G)$ is obviously monotone in the sense that $F\subseteq G$ implies $C_{\rm OR}(F)\le C_{\rm OR}(G)$ and the graphs $M_r$ are triangle-free themselves, we have that
$\sup\{C_{\rm OR}(G): K_3\nsubseteq G\}=\lim_{r\to\infty}C_{\rm
  OR}(M_r)$
and thus the statement of the theorem. (Note that the existence of the limit
follows from $M_{r-1}\subseteq M_r$ and thus $C_{\rm OR}(M_r)\ge C_{\rm OR}(M_{r-1})$ being true for every $r$.)
\hfill$\Box$

\smallskip
\par\noindent
We do not know how to show that the Shannon OR-capacity of Mycielski graphs does or does not go to infinity. As a modest step forward we give a lower bound on the Shannon OR-capacity of the Gr\"otzsch graph that is strictly larger than $C_{\rm OR}(C_5)=\sqrt{5}$, i.e, the Shannon OR-capacity of the previous Mycielski graph.

\begin{thm}
$$C_{\rm OR}(M_4)\ge \sqrt[4]{28}.$$
\end{thm}

\proof
As $C_{\rm OR}(G)\ge\sqrt[t]{\omega(G^t)}$ for any finite $t$, it is enough to present a clique of size $28$ in the $4$th OR-power of the Gr\"otzsch graph $M_4=M(C_5)$. This is what we do.
We denote the vertices of $C_5$ by $0,1,2,3,4$, where every vertex is connected to the vertices labeled by neighboring integers, plus $\{0,4\}$ is also an edge. The corresponding twin vertices (cf. Definition~\ref{defi:Myc}) are labeled $0',1',2',3',4'$, respectively, and the eleventh vertex of $M(C_5)$ is labeled $z$ as in Definition~\ref{defi:Myc}.

\bigskip
\par\noindent
We claim that the following $28$ vertices form a clique in $[M(C_5)]^2$.
Indeed, it is enough to observe, that the $25$ vertices listed in the first
five rows below are created from Shannon's well-known $5$-element clique of
$C_5^2$ in \cite{Sha} given as $00, 12,24, 31, 43$ by taking its second power
and then "lifting" some vertices to their twins in such a way that at least
one connection still remained between any two of our sequences. All sequences
except five have a lifted element in the last two coordinates and for the five
exceptional ones (see them put diagonally) both of the first two coordinates are lifted. Thus all these sequences will be connected with the three in the last row that are also connected to each other.

$$0'0'00\ \ \ 120'0\ \ \ 2400'\ \ \ 3100'\ \ \ 430'0$$
$$001'2\ \ \ 1'2'12\ \ \ 241'2\ \ \ 3112'\ \ \ 4312'$$
$$0024'\ \ \ 122'4\ \ \ 2'4'24\ \ \ 312'4\ \ \ 4324'$$
$$0031'\ \ \ 1231'\ \ \ 243'1\ \ \ 3'1'31\ \ \ 433'1$$
$$004'3\ \ \ 1243'\ \ \ 2443'\ \ \ 314'3\ \ \ 4'3'43$$
$$0'zzz\ \ \ \ z0'zz\ \ \ \ zzzz$$

\smallskip
\par\noindent
Thus we have $\omega([M(C_5)]^2)\ge 28$ and the statement follows.
\hfill$\Box$


\begin{thebibliography}{99}



\bibitem{AKM} Noga Alon, J\'anos K\"orner, Angelo Monti, String quartets in binary, {\em Combin., Probab., Comput.}, 9 (2000), 381--390.  

\bibitem{AO}
Noga Alon and Alon Orlitsky, Repeated communication and Ramsey graphs, {\em IEEE
Trans. Inform. Theory}, 41 (1995), 1276-1289.

\bibitem{AxFurMub} Maria Axenovich, Zolt\'an F\"uredi, Dhruv Mubayi, On
  generalized Ramsey theory: the bipartite case, {\em J. Combin. Theory
    Ser. B}, 79 (2000), no. 1, 66--86.

\bibitem{AxIv} Maria Axenovich and Perry Iverson, Edge-colorings avoiding rainbow
  and monochromatic subgraphs, {\em Discrete Math.}, 308 (2008), no. 20,
  4710--4723.

\bibitem{BFGyR} Aart Blokhuis, Ralph Faudree, Andr\'as Gy\'arf\'as, Mikl\'os
  Ruszink\'o, Anti-Ramsey colorings in several rounds, {\em J. Combin. Theory
    Ser. B}, 82 (2001), no. 1, 1--18.

\bibitem{BH} Tom Bohman and Ron Holzman, A nontrivial lower bound on the Shannon capacities of the complements of odd cycles, {\em IEEE Trans. Inform. Theory}, 49 (2003), 721--722.

\bibitem{Bondy} J.\ A. Bondy, Basic graph theory: paths and circuits, in: {\em
  Handbook of combinatorics, Vol. 1}, (Ronald R. Graham, Martin Gr\"otschel,
  L\'aszl\'o Lov\'asz eds.) 3--110, Elsevier, Amsterdam, 1995.

\bibitem{CG} Fan Chung and Ron Graham, {\em Erd\H{o}s on Graphs. His Legacy of Unsolved Problems},
A K Peters/CRC Press, New York, 1998.

\bibitem{CGyL} Mathew Cropper, Andr\'as Gy\'arf\'as, Jen\H{o} Lehel, Hall
  ratio of the Mycielski graphs, {\em Discrete Math.}, 306 (2006), no. 16,
  1988--1990.

\bibitem{Radha} Marco Dalai, Venkatesan Guruswami, Jaikumar Radhakrishnan, An improved bound on
the zero-error list-decoding capacity of the 4/3 channel, {\em Proc. 2017 IEEE
International Symposium on Information Theory, ISIT 2017}, 1658--1662,
https://doi.org/10.1109/ISIT.2017.8006811

\bibitem{PE} Peter Elias, Zero-error capacity under list decoding, {\em IEEE Trans. Inform. Theory,}, 34 (1988), 1070--1074.

\bibitem{Erdos} Paul Erd\H{o}s, Solved and unsolved problems in combinatorics
  and combinatorial number theory, Proceedings of the Twelfth Southeastern
  Conference on Combinatorics, Graph Theory and Computing, Vol. I (Baton
  Rouge, La., 1981). {\em Congr. Numer.}, 32 (1981), 49--62.

\bibitem{ErdGyar} Paul Erd\H{o}s and Andr\'as Gy\'arf\'as, A variant of the
  classical Ramsey problem, {\em Combinatorica}, 17 (1997), no. 4, 459--467.

\bibitem{EMcT} Paul Erd\H{o}s, Robert J. McEliece, Herbert Taylor, Ramsey bounds for graph products. {\em Pacific J. Math.}, 37 (1971), 45--46.

\bibitem{Exoo} Geoffrey Exoo, A lower bound for Schur numbers and multicolor Ramsey numbers of $K_3$, {\em Electronic J. Combin.}, 1 (1994), \# R8.

\bibitem{FK} Michael L. Fredman and J\'anos Koml\'os, On the size of
  separating systems and families of perfect hash functions, {\em SIAM J. Algebraic Discrete Methods}, 5 (1984), no. 1, 61--68.

\bibitem{GKV} Luisa Gargano, J\'anos K\"orner, Ugo Vaccaro, Qualitative independence and Sperner problems for directed graphs, {\em J. Combin. Theory Ser A}, 61 (1992), 173--192.

\bibitem{HRCW} Haim Hanani, Dijen K. Ray-Chaudhury, Richard M. Wilson, On
  resolvable designs, {\em Discrete Math.}, 3 (1972), 343--357.

\bibitem{HS} Zita Helle and G\'abor Simonyi, Orientations making k-cycles
  cyclic, {\em Graphs Combin.}, 32 (2016), 2415--2423.

\bibitem{KHGP} J\'anos K\"orner, Fredman-Koml\'os bounds and information theory,
{\em SIAM J. Algebraic Discrete Methods}, 7 (1986), no. 4, 560--570.

\bibitem{KM} J\'anos K\"orner and Katalin Marton, New bounds for perfect
  hashing via information theory, {\em European J. Combin.}, 9 (1988), no. 6,
  523--530.

\bibitem{KS1} J\'anos K\"orner and G\'abor Simonyi, Separating partition systems and very different sequences, {\em SIAM J. Discrete Mathematics}, 1 (1988), 355-359.

\bibitem{KS} J\'anos K\"orner and G\'abor Simonyi, Trifference, {\em Studia Sci. Math. Hungar.}, 30 (1995), 95--103, and also in: ''{\em Combinatorics and its Applications to the Regularity and Irregularity of Structures}'', Walter A. Deuber and Vera T. S\'os eds., Akad\'emiai Kiad\'o, Budapest, 1995.

\bibitem{LovKne}  L\'aszl\'o Lov\'asz, Kneser's conjecture, chromatic number, and
homotopy, {\em J. Combin.\ Theory Ser.\ A}, 25 (1978), no.\ 3,
319--324.

\bibitem{LovShan} L\'aszl\'o Lov\'asz, On the Shannon capacity of graphs, {\em IEEE Trans. Inform. Theory}, 25 (1979), 1--7.

\bibitem{Mubayi} Dhruv Mubayi, An explicit construction for a Ramsey
  problem, {\em Combinatorica}, 24 (2004), no. 2, 313--324.

\bibitem{Myc} Jan Mycielski, Sur le coloriage des graphs. (French) {\em Colloq. Math.} 3, (1955). 161--162.

\bibitem{NR} Jaroslav Ne\v{s}et\v{r}il and Moshe Rosenfeld, I. Schur, C. E. Shannon and Ramsey
numbers, a short story. Combinatorics, graph theory, algorithms and
applications. {\em Discrete Math.} 229 (2001), no. 1-3, 185--195.

\bibitem{RCW} Dijen K. Ray-Chaudhury and Richard M. Wilson, The existence of
  resolvable block designs,{\em Survey of combinatorial theory
    (Proc. Internat. Sympos., Colorado State Univ., Fort Collins, Colo.,
    1971)}, North-Holland, Amsterdam, 1973, 361--375.

\bibitem{SaSi} Attila Sali and G\'abor Simonyi, Orientations of self-complementary graphs
  and the relation of Sperner and Shannon capacities, {\em Europ.
    J. Combin.}, {\bf 20} (1999), 93--99.

\bibitem{Sha} Claude E. Shannon, The zero-error capacity of a noisy channel, {\em IRE Trans. Inform. Theory}, 2 (1956), 8--19.

\bibitem{Stinson} Douglas R. Stinson, A survey of Kirkman triple systems and
  related designs, {\em Discrete Math.}, 92 (1991), no. 1-3, 371--393.


\end{thebibliography}
\end{document}